\documentclass[11pt,twoside]{article}
\usepackage[hmargin=1.in, vmargin=1.in, a4paper, centering]{geometry}
\usepackage[T1]{fontenc}
\usepackage[lining]{ebgaramond}
\usepackage{amsmath,amsthm,mathtools}
\usepackage[ebgaramond]{newtxmath}
\usepackage{enumerate}
\usepackage[font=footnotesize, labelfont=bf, labelsep=colon, margin=0.5in]{caption}
\usepackage{graphicx}
\graphicspath{{figures/}}

\usepackage{xcolor}
\usepackage[toc,titletoc]{appendix}
\appendixtocoff
\usepackage{etoolbox, apptools}
\AtBeginEnvironment{appendices}{\appendixtrue}
\usepackage[nobottomtitles,pagestyles]{titlesec}
\titleformat{\section}
{\normalfont\large\bfseries}
{\filcenter\IfAppendix{Appendix }{\S}\thesection.}{1ex}{\filcenter}
\titleformat{\subsection}
{\normalfont\bfseries}
{\filcenter\S\thesubsection.}{1ex}{\filcenter}
\usepackage[colorlinks,allcolors=blue,pagebackref]{hyperref}
\usepackage{pifont}
\renewcommand*{\backrefalt}[4]{%
\ifcase #1 %
\textcolor{red}{No citations}%
\or
\ding{43}~p.~#2%
\else
\ding{43}~pp.~#2%
\fi}
\usepackage{breakcites}
\newcommand{\mydoi}[1]{\href{https://doi.org/#1}{doi: #1}}
\newcommand{\myarXiv}[1]{\href{https://arxiv.org/abs/#1}{arXiv: #1}}
\numberwithin{equation}{section}
\ifdefined\iflatexml
\theoremstyle{plain}
\newtheorem{theorem}{Theorem}[section]
\newtheorem{proposition}[theorem]{Proposition}
\newtheorem{lemma}[theorem]{Lemma}
\newtheorem{corollary}[theorem]{Corollary}
\theoremstyle{definition}

\newtheorem{remark}[theorem]{Remark}
\else
\usepackage{tcolorbox,varwidth}
\tcbuselibrary{skins, breakable, theorems}
\colorlet{thmcolor}{orange!80!white}
\colorlet{conjcolor}{orange!40!white}
\colorlet{remcolor}{teal!50!white}
\colorlet{defncolor}{blue!50!black}
\usepackage{keytheorems}
\tcbset{commonstyle/.style={%
colback=white!95!tcbcolframe,
colbacktitle=white!80!tcbcolframe,
arc=2mm,
fonttitle=\bfseries,
coltitle=black,
enhanced,
varwidth boxed title*=-2cm,
attach boxed title to top left={yshift=-3mm,xshift=0.5cm,yshifttext=-1mm},
description font=\mdseries,
breakable,
before upper={\parindent15pt\noindent},
beforeafter skip balanced=10.0pt plus 1.0pt minus 1.0pt,
separator sign={\textbf{: }}%,
overlay first={\draw[color=tcbcolframe, line width=.5pt] (frame.south west)--(frame.south east);},
overlay middle={\draw[color=tcbcolframe, line width=.5pt] (frame.south west)--(frame.south east);\draw[color=tcbcolframe, line width=.5pt] (frame.north west)--(frame.north east);},
overlay last={\draw[color=tcbcolframe, line width=.5pt] (frame.north west)--(frame.north east);}
}
}
\newkeytheoremstyle{mythmstyle}{
bodyfont=\normalfont,
headpunct={},
notebraces={}{},
noteseparator={: }, 
notefont=\bfseries,
tcolorbox = {colframe=thmcolor,commonstyle}
}
\newkeytheoremstyle{myremstyle}{
bodyfont=\normalfont,
headpunct={},
notebraces={}{},
noteseparator={: }, 
notefont=\bfseries,
tcolorbox = {colframe=remcolor,commonstyle}
}
\newkeytheoremstyle{mydefstyle}{
bodyfont=\normalfont,
headpunct={},
notebraces={}{},
noteseparator={: }, 
notefont=\bfseries,
tcolorbox = {colframe=defncolor,commonstyle}
}
\newkeytheorem{theorem}[
name=Theorem,
style=mythmstyle,
parent=section
]
\newkeytheorem{lemma}[
name=Lemma,
style=mythmstyle,
sibling=theorem
]
\newkeytheorem{proposition}[
name=Proposition,
style=mythmstyle,
sibling=theorem
]
\newkeytheorem{corollary}[
name=Corollary,
style=mythmstyle,
sibling=theorem
]
\newkeytheorem{remark}[
name=Remark,
sibling=theorem,
style=myremstyle
]
\newkeytheorem{definition}[
name=Definition,
style=mydefstyle,
sibling=theorem
]
\fi
\usepackage{mleftright}
\mleftright
\newcommand{\R}{\mathbb R}
\newcommand{\C}{\mathbb C}
\newcommand{\D}{\mathbb D}
\newcommand{\eps}{\varepsilon}
\newcommand{\usc}{u^{\mathrm{sc}}}
\newcommand{\Pcal}{\mathcal P}
\newcommand{\Hcal}{\mathcal H}
\newcommand{\Fcal}{\mathcal F}
\newcommand{\Ycal}{\mathcal Y}
\newcommand{\Xcal}{\mathcal X}
\newcommand{\id}{\mathrm{id}}
\newcommand{\Ran}{\operatorname{Ran}}
\newcommand{\ir}{\mathrm{i}}
\newcommand{\dr}{\mathrm{d}}
\renewcommand{\tilde}{\widetilde}
\allowdisplaybreaks[2]
\newcommand{\asize}{\footnotesize}
\title{A family of homogeneous convex non-scattering domains with analytic non-circular  boundaries}
\author{Fioralba Cakoni\\%
\asize Department of Mathematics, Rutgers University\\%
\asize New Brunswick, NJ 08901-8554, USA\\%
\asize \href{mailto:fc292@math.rutgers.edu}{\nolinkurl{ fc292@math.rutgers.edu}}\\%
\asize \url{https://sites.math.rutgers.edu/~fc292}\\%
\asize ORCID: 0000-0002-4648-8780
\and 
Michael Levitin\\
\asize Department of Mathematics and Statistics,  University of Reading\\%
\asize Pepper Lane, Whiteknights, Reading RG6 6AX, U. K.\\%
\asize \href{mailto:M.Levitin@reading.ac.uk}{\nolinkurl{M.Levitin@reading.ac.uk}}\\%
\asize \url{https://www.michaellevitin.net}\\%
\asize ORCID: 0000-0003-0020-3265}
\date{\footnotesize 19 November 2026}
\begin{document}
\maketitle
\thispagestyle{empty}
\begin{abstract}
A homogeneous planar  non-scattering domain  is a bounded  medium with a constant material coefficient that produces no scattered field for some nonzero entire solution of the Helmholtz equation at a particular wave number. It is known that the boundary of such a domain must be real analytic at all points where the incident field does not vanish. Previously, disks were the only known examples. We construct a family of homogeneous strictly convex non-scattering domains that are neither disks nor ellipses. The family bifurcates from a particular non-scattering disk, with the boundary, wave number, and constant coefficient varying along the branch. The construction uses a conformal change of variables and analytic bifurcation theory applied to the resulting nonlinear problem.\end{abstract}

\begin{center}
{\footnotesize\noindent

}
\end{center}

{\footnotesize\tableofcontents}

%%%%%%%%%%%%%%%%%%%%%%%%%%%%%%%%%%%%%%%%%%%%%%%%%%%%%%%%%%%%%%%%
\section{Formulation of the problem and main results}\label{sec:model}

\subsection{Non-scattering domains}

The question investigated in this paper arises from the scattering of
time-harmonic waves by a homogeneous penetrable medium of bounded support,
modelled by the reduced wave equation. More precisely, let
$\Omega\subset\R^2$ denote the bounded domain occupied by the medium,
with a constant positive coefficient $\sigma$ (the square of the refractive index), $\sigma\ne1$; the refractive
index outside $\Omega$ is normalised to one. If $k>0$ is the exterior
wave number, which is proportional to the probing frequency, then the
interior wave number is
\begin{equation}\label{eq:interior}
\kappa=k\sqrt{\sigma},\qquad\text{so that}\qquad \sigma=\frac{\kappa^2}{k^2}.
\end{equation}
Let $v$ be an entire incident field satisfying
\begin{equation}\label{eq:incident-equation}
(\Delta+k^2)v=0\qquad\text{in }\R^2.
\end{equation}
The corresponding  interior total field $u$ and the exterior scattered field $\usc$
solve the scattering problem
\begin{equation}\label{eq:scattering}
\begin{cases}
(\Delta+\kappa^2)u=0&\text{in }\Omega,\\
(\Delta+k^2)\usc=0&\text{in }\R^2\setminus\overline\Omega,\\
u=v+\usc&\text{on }\partial\Omega,\\
\partial_\nu u=\partial_\nu v+\partial_\nu \usc
&\text{on }\partial\Omega,\\
\displaystyle\lim_{r\to\infty}\sqrt r\,
\left(\partial_r \usc-\ir k \usc\right)=0,
\end{cases}
\end{equation}
where $\nu$ is the outward unit normal to $\partial\Omega$, $r=|x|$, and the
Sommerfeld radiation condition in the last line holds uniformly with
respect to $x/|x|\in\mathbb S^1$. We call $\Omega$ a
\emph{non-scattering domain} at wave number $k$ for the coefficient
$\sigma$ if $\usc\equiv0$ for some nonzero entire incident field $v$.
Equivalently, there is a pair $(u,v)$, with $v$ a nonzero entire solution
of \eqref{eq:incident-equation}, satisfying the following homogeneous problem
\begin{equation}\label{eq:transmission}
\begin{cases}
(\Delta+\kappa^2)u=0&\text{in }\Omega,\\
(\Delta+k^2)v=0&\text{in }\Omega,\\
u=v&\text{on }\partial\Omega,\\
\partial_\nu u=\partial_\nu v&\text{on }\partial\Omega.
\end{cases}
\end{equation}

The unit disk $\D=\{x\in\R^2:|x|<1\}$ is non-scattering for suitable
$k\ne\kappa$. Indeed, the radial fields
\begin{equation}\label{eq:disk-fields}
v(x)=J_0(k|x|),\qquad
u(x)=\frac{J_0(k)}{J_0(\kappa)}J_0(\kappa|x|)
\end{equation}
satisfy \eqref{eq:transmission} whenever $J_0(k)J_0(\kappa)\ne0$ and
\begin{equation}\label{eq:diskmatch}
d_0(k)=d_0(\kappa),\qquad
d_0(t):=\frac{tJ_0'(t)}{J_0(t)}.
\end{equation}
For a fixed coefficient $\sigma>0$, $\sigma\ne1$, the non-scattering
wave numbers $k>0$ satisfying \eqref{eq:diskmatch}, with
$\kappa=k\sqrt{\sigma}$, form a countably infinite set with $+\infty$
as the only accumulation point \cite[\S8.6]{ColtonKress2019}.
Translation and dilation give corresponding non-scattering
configurations for any disk. Similar conclusions hold for radially
varying refractive indices.

Non-scattering is closely connected with the transmission eigenvalue
problem, which is the homogeneous system \eqref{eq:transmission} posed for $u\in L^2(\Omega)$ and $v\in L^2(\Omega)$
with $u-v\in H^2_0(\Omega)$. Here $k$ is a transmission eigenvalue if
\eqref{eq:transmission} has a non-trivial solution pair $(u,v)$ defined
in $\Omega$, called a transmission eigenfunction pair. If an incident
field $v$, which is an entire solution of
\eqref{eq:incident-equation}, produces no scattered field outside
$\Omega$, then its restriction $v|_{\Omega}$ and the transmitted field
$u$ form a transmission eigenfunction pair. Conversely, for a transmission eigenfunction pair to give a non-scattering
configuration, its background component $v$ must
extend to an admissible incident field satisfying the Helmholtz equation
in the whole space. Of particular interest in scattering theory are
Herglotz incident fields, which are superpositions of plane waves.
In this case, the extension must be a Herglotz wave function.
It is well known that a
homogeneous medium $(\Omega,\sigma)$, $\sigma\neq 1$, has infinitely many
transmission eigenvalues $k>0$, with $+\infty$ as the only accumulation
point, irrespective of the regularity of $\Omega$.
Thus, the existence of real transmission eigenvalues does not by itself
imply non-scattering, although the non-scattering wave numbers form a
(possibly empty) subset of the transmission eigenvalues and hence of a discrete set. We refer to
\cite[Chapter~7]{CCH} or the survey paper \cite{CakoniVogelius2026}
for a systematic account of this connection.

For a spatially constant coefficient $\sigma>0$, $\sigma\ne1$,
non-scattering imposes strong restrictions on the boundary of the
inhomogeneity. The free-boundary regularity results in \cite{CV,SS}
imply that a Lipschitz boundary must be real analytic near every
boundary point where the incident field does not vanish. In particular,
if the incident field is nowhere zero on the boundary, then the entire
boundary is real analytic.

Scattering by analytic inhomogeneities with constant material
coefficients has been investigated in several works~\cite{HV1,HV2,HV3}.
For suitable classes of analytic boundaries, a complex-analytic approach
establishes scattering whenever the analytically continued incident
field does not vanish at an accessible critical point determined by the
complexified boundary~\cite{HV1}. Subsequent work establishes
plane-wave scattering criteria involving the frequency, incident
direction, and geometry~\cite{HV2}. Most recently, it was proved that
every incident plane wave scatters from a bounded, strictly convex
planar inhomogeneity with constant material parameters~\cite{HV3}.
An infinite hierarchy of necessary differential conditions for
non-scattering at points determined by the complexified geometry was
further derived in \cite{H}.

Complementary finiteness results were obtained in \cite{VX21,VX25}.
For classes of non-circular perturbations of disks satisfying the
geometric hypotheses in those papers, Herglotz incident fields of the form
\begin{equation}\label{eq:incidentG}
v_{g,k}(x)=
\frac1{2\pi}\int_0^{2\pi} g(\varphi)
e^{\ir kx\cdot(\cos\varphi,\sin\varphi)}\,\dr\varphi,
\end{equation}
with a fixed, nonzero, sufficiently regular Herglotz density $g$, can be
non-scattering for at most finitely many positive wave numbers.
These results do not exclude individual non-scattering configurations.

In this paper we construct a family of small deformations of the unit
disk that are strictly convex, have real-analytic boundaries, and are
neither disks nor ellipses. Constructing non-scattering homogeneous media, that is, media with a
spatially constant material coefficient, is a challenging problem. To our knowledge, the domains constructed in the present paper provide
the first non-circular, simply connected examples of non-scattering
homogeneous media with $\sigma\neq1$.

Each domain we construct is non-scattering for a radial Herglotz
wave of the form $v(x)=J_0(k|x|)$. The shape of the domain, the wave
numbers, and the coefficient $\sigma$ vary together along the family.
Our construction starts from a non-scattering disk for which $\kappa^2$ is a particular
Dirichlet eigenvalue of the disk. This resonance permits an infinitesimal
non-circular deformation, which bifurcation theory extends to a family
of exact nonlinear solutions parametrised by conformal maps of the disk. Exact integral identities then determine the leading corrections to the polar radius
and material coefficient and show that the domains on the branch are not ellipses.

Setting $w:=u-v$ in \eqref{eq:transmission}, we see that a homogeneous medium $(\Omega,\sigma)$, with $\sigma\neq1$, is non-scattering if and only if the problem
\begin{equation}\label{eq:nonscat}
\begin{cases}
\Delta w+\kappa^2 w=(k^2-\kappa^2)v&\text{in }\Omega,\\
w=0&\text{on }\partial\Omega,\\
\partial_\nu w=0&\text{on }\partial\Omega
\end{cases}
\end{equation}
has a solution for some $k>0$ and some non-trivial entire incident field $v$ satisfying $(\Delta+k^2)v=0$ (recall that $\kappa=k\sqrt{\sigma}$).

\begin{remark}
There is a close parallel with the study of \emph{Schiffer domains}, or equivalently, for smooth 
simply connected domains, those domains which lack the \emph{Pompeiu property}. More specifically, a
domain $\Omega$ is a Schiffer domain if  there exist $\mu>0$ and a non-constant Neumann Laplacian eigenfunction $\phi$ which is constant on the boundary,
\[
\begin{cases}
\Delta \phi +\mu \phi=0&\text{in }\Omega,\\
\phi=1&\text{on }\partial\Omega,\\
\partial_\nu \phi=0&\text{on }\partial\Omega.
\end{cases}
\]
Setting $\mu:=\kappa^2$ and $w:=\phi-1$, we get
\begin{equation}\label{schieffer}
\begin{cases}
\Delta w + \kappa^2 w= -\kappa^2 &\text{in }\Omega,\\
w=0&\text{on }\partial\Omega,\\
\partial_\nu w=0&\text{on }\partial\Omega.
\end{cases}
\end{equation}

The Schiffer conjecture stated that disks are the only bounded simply connected
planar domains with smooth boundaries for which  \eqref{schieffer} holds. This conjecture was recently disproved, independently
and almost simultaneously, in \cite{CaoLaboraDeDiosPont2026} and
\cite{ColbrookStepaniants2026}, where non-circular Schiffer domains were constructed.

The comparison of \eqref{schieffer} and \eqref{eq:nonscat} immediately shows that any Schiffer domain is ``non-scattering'' with $k=0$, or, equivalently, for a constant incident field.  As our definition of non-scattering domains requires $k>0$, the constructed Schiffer domains are not automatically non-scattering in this sense.
\end{remark}

\subsection{Main results}
Throughout, the incident field is the radial Herglotz wave function \eqref{eq:incidentG} with constant
density $g\equiv 1$,
\begin{equation}\label{eq:incident}
v_k(x):=v_{1,k}(x)=J_0(k|x|)
=\frac1{2\pi}\int_0^{2\pi}
e^{\ir kx\cdot(\cos\varphi,\sin\varphi)}\,\dr\varphi,
\end{equation}
a superposition of plane waves from all directions. Taking $v=v_k$ in \eqref{eq:nonscat}, we obtain the free-boundary formulation of non-scattering in our setting:
\begin{equation}\label{eq:clamped}
\begin{cases}
(\Delta+\kappa^2)w=(k^2-\kappa^2)J_0(k|x|)&\text{in }\Omega,\\
w=0&\text{on }\partial\Omega,\\
\partial_\nu w=0&\text{on }\partial\Omega.
\end{cases}
\end{equation}
This is a free-boundary problem: we seek a domain $\Omega$ on whose
boundary both Cauchy traces of a solution $w\in H^2(\Omega)$ to the
inhomogeneous source problem vanish. If $w\in H^2(\Omega)$
solves \eqref{eq:clamped}, then $(u,v)=(v_k+w,v_k)$ solves
\eqref{eq:transmission}, and $\Omega$ is non-scattering for the
incident field $v_k$.

Let $j_{n,1}$ denote the first positive zero of $J_n$. The disk from which
we bifurcate has the parameters
\begin{equation}\label{eq:seeds}
\kappa_*=j_{2,1},\qquad
d_0(k_*)=d_0(\kappa_*),\qquad 0<k_*<j_{0,1},
\end{equation}
which are well defined, see Lemma~\ref{lem1}. Numerically,
\begin{align}
\kappa_*&=5.1356223018406826\ldots,\nonumber\\
k_*&=2.2308555506296499\ldots,\label{eq:seed-values}\\
\sigma_*:=\kappa_*^2/k_*^2&=5.2996019548130573\ldots.\nonumber
\end{align}

Throughout the paper, we fix a H\"older exponent $\alpha\in(0,1)$.

\begin{theorem}\label{thm:non-scattering-branch}
There exist $\eps_0>0$ and three real-analytic curves
on $I_0=(-\eps_0,\eps_0)$,
\[
\begin{aligned}
I_0\ni\eps&\longmapsto
\Phi_\eps\in C^{4,\alpha}(\overline{\D};\C),
&\qquad \Phi_0&=\id,\\
I_0\ni\eps&\longmapsto k_\eps\in(0,\infty),
& k_0&=k_*,\\
I_0\ni\eps&\longmapsto\kappa_\eps\in(0,\infty),
& \kappa_0&=\kappa_*,
\end{aligned}
\]
(where $\C$ is identified with $\R^2$) with the following properties:
\begin{itemize}
\item For every $\eps\in I_0$, the map $\Phi_\eps$
is holomorphic in $\D$, injective on
$\overline{\D}$, and normalised by
$\Phi_\eps(0)=0$ and $\Phi_\eps'(0)>0$. The domain $\Omega_\eps=\Phi_\eps(\D)$
has real-analytic boundary with everywhere positive curvature
and is invariant under reflection in both coordinate axes.
\item With $\sigma_\eps=\kappa_\eps^2/k_\eps^2>1$,
there exists a real-valued function $u_\eps\in H^2(\Omega_\eps)$ such that
$(u_\eps,v_{k_\eps})$ satisfies
\eqref{eq:transmission} on $\Omega_\eps$ with wave
numbers $k_\eps$ and $\kappa_\eps$.
\end{itemize}
Thus $\Omega_\eps$ with spatially constant coefficient
$\sigma_\eps$ is non-scattering for the  Herglotz
incident wave $v_{k_\eps}$ defined by \eqref{eq:incident}.
For every $\eps\in I_0\setminus\{0\}$, the domain
$\Omega_\eps$ is neither a disk nor an ellipse.
\end{theorem}

\begin{figure}[htb]
\centering
\includegraphics{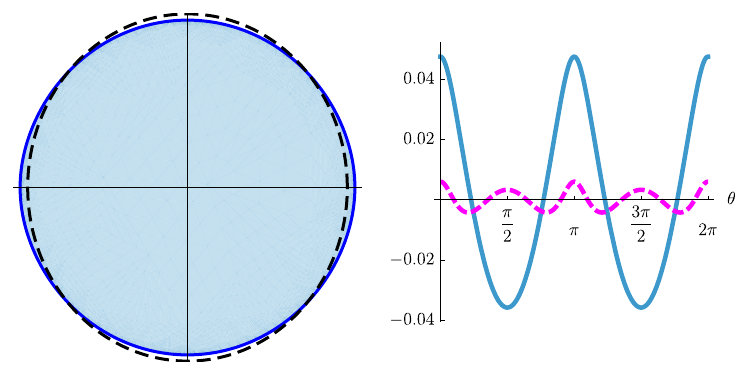}
\caption{Numerical approximation to the non-scattering domain at
$\eps=0.04$. The boundary was computed from the truncated Fourier series
$\rho_\eps(\theta)=1+\eps\cos2\theta+\sum_{j=2}^{24}\hat\rho_j\cos2j\theta$,
with $\hat\rho_j$, $k_\eps$ and $\kappa_\eps$ determined numerically from the
moment identities \eqref{eq:radial-moments} for $n=0,2,\dots,48$.
Left: the computed boundary (solid) and the unit circle (dashed).
Right: the radial differences $\rho_\eps(\theta)-1$ (solid) and
$\rho_\eps(\theta)-\rho_E(\theta)$ (dashed), plotted on the same scale.
Here $\rho_E$ is the polar radius of the centred, axis-aligned ellipse with semi-axes $1.04151$ and $0.960921$, 
obtained by minimising the mean-square difference of polar radii over
$\theta\in[0,2\pi]$.}
\label{fig:nonscattering-domain}
\end{figure}

\begin{remark}
Numerical experiments indicate that one can take $\eps_0\approx0.0505$ in Theorem~\ref{thm:non-scattering-branch}; at this value, $\partial\Omega_\eps$
reaches the first nodal circle of the incident field $v_{k_\eps}$, and the boundary
appears to lose smoothness at the two points of contact. For a numerical approximation of $\Omega_\eps$ with $\eps=0.04$ see Figure~\ref{fig:nonscattering-domain}. We also note that, by analytic elliptic regularity, $u_\eps$ is real analytic up to $\partial\Omega_\eps$. 
\end{remark}

\begin{proposition}\label{prop:normalised-expansions}
The family constructed in the proof of Theorem~\ref{thm:non-scattering-branch} satisfies
$$\Omega_\eps
=\left\{r(\cos\theta,\sin\theta):
0\le r<\rho_\eps(\theta)\right\},$$
where the polar radius is normalised by
\begin{equation}\label{eq:polar-normalisation}
\frac1{2\pi}\int_0^{2\pi}\rho_\eps(\theta)\,\dr\theta=1,
\qquad
\frac1\pi\int_0^{2\pi}\rho_\eps(\theta)\cos2\theta\,\dr\theta
=\eps.
\end{equation}
As $\eps\to0$,
\begin{equation}\label{eq:mainshape}
\rho_\eps(\theta)=1+\eps\cos2\theta+
\left(\frac34+\frac{\kappa_*^2}{12}\right)\eps^2\cos4\theta
+O_{C^0(\mathbb S^1)}(\eps^3),
\end{equation}
and the material coefficient satisfies
\begin{equation}\label{eq:main-contrast}
\begin{aligned}
\sigma_\eps
&=\sigma_*+c_*\eps^2+O(\eps^4),\\
c_*&=
\frac{2\sigma_*(2\kappa_*^2+3)(\kappa_*^2-k_*^2)}
{3(\kappa_*^4+4k_*^2)}>0.
\end{aligned}
\end{equation}
\end{proposition}

To first order, $\Omega_\eps$ is stretched along one coordinate axis and
compressed along the other, and so resembles an ellipse. The difference
appears at second order: for an ellipse normalised as in
\eqref{eq:polar-normalisation}, the coefficient of $\eps^2\cos4\theta$ is
$3/4$ rather than $3/4+\kappa_*^2/12$.

The following corollary of Theorem~\ref{thm:non-scattering-branch} and
Proposition~\ref{prop:normalised-expansions} shows that both the coefficient
$\sigma$, in a suitable interval, and the exterior wave number $k>0$ may be
prescribed.
\begin{corollary}
\label{cor:prescribed-parameters}
There exists $\delta>0$ such that, for every
$\sigma\in(\sigma_*,\sigma_*+\delta)$ and every $k>0$,
there is a bounded, simply connected domain
$\Omega\subset\R^2$ containing the origin,
with real-analytic boundary of everywhere positive
curvature, invariant under reflection in both coordinate
axes, and neither a disk nor an ellipse, for which
there exists a real-valued $u\in H^2(\Omega)$ such that
$(u,v_k)$ satisfies \eqref{eq:transmission} with
$\kappa=k\sqrt\sigma$. Thus $\Omega$ with spatially
constant coefficient $\sigma$ is non-scattering for the
Herglotz incident wave \eqref{eq:incident}.
\end{corollary}
\begin{proof}
By Theorem~\ref{thm:non-scattering-branch} and Proposition~\ref{prop:normalised-expansions}, $\sigma_\eps$ is real analytic on $(-\eps_0,\eps_0)$ and satisfies \eqref{eq:main-contrast}, so that
\[
\frac{\dr\sigma_\eps}{\dr\eps}=2c_*\eps+O(\eps^3),\qquad c_*>0.
\]
Hence there exists $\eps_1\in(0,\eps_0)$ such that $\sigma_\eps$ is strictly increasing on $[0,\eps_1]$. Set $\delta=\sigma_{\eps_1}-\sigma_*>0$. Given
$\sigma\in(\sigma_*,\sigma_*+\delta)$, there is a unique $\eps\in(0,\eps_1)$ such that $\sigma_\eps=\sigma$.
The corresponding domain $\Omega_\eps$ is non-scattering at the exterior wave number $k_\eps$. For the prescribed exterior wave number $k>0$, set
\[
\lambda=\frac{k_\eps}{k},
\qquad
\Omega=\lambda\Omega_\eps,
\qquad
u(x)=u_\eps(x/\lambda).
\]
The interior wave number on the dilated domain is
\[
\kappa=\frac{\kappa_\eps}{\lambda}
=k\frac{\kappa_\eps}{k_\eps}
=k\sqrt{\sigma}.
\]
With $y=x/\lambda$, we have
$v_k(x)=J_0(k_\eps|y|)=v_{k_\eps}(y)$.
The dilation scales the Laplacian by $\lambda^{-2}$ and both
normal derivatives by $\lambda^{-1}$, so the Helmholtz equations
hold at the rescaled wave numbers and the Cauchy data match on
$\partial\Omega$. Thus $u$ is real-valued, belongs to
$H^2(\Omega)$, and $(u,v_k)$ satisfies \eqref{eq:transmission}.
Dilation preserves all the stated geometric properties. For a fixed
coefficient $\sigma$, the domains obtained by varying $k$ are dilates of
the same normalised domain $\Omega_\eps$, with mean polar
radius $\lambda=k_\eps/k$.
\end{proof}

\begin{remark}
The domains $\Omega_\eps$ are exactly non-scattering for the associated
wave numbers and material coefficients. The expansion \eqref{eq:mainshape} describes
these exact domains and is used only to establish non-ellipticity; exact
non-scattering is not asserted for the boundary obtained by discarding
the remainder.
\end{remark}

\begin{remark}
The scattering criterion in \cite{HV1} requires both geometric assumptions and a non-vanishing condition on the analytically continued incident field at a suitable complex critical point. At least one of these assumptions must fail for the non-scattering configurations constructed here; determining which one fails remains open. Positivity of the incident field on the real boundary does not imply its non-vanishing at the complex critical point.
\end{remark}

The paper is organised as follows.  In \S\ref{main} we construct the unnormalised branch of
non-scattering domains. We first choose the resonant interior wave
number $\kappa_*=j_{2,1}$ and use the radial matching condition to
obtain an analytic curve $k=k(\kappa)$ along which the unit disk is
non-scattering (\S\ref{sec:disk}). This gives the trivial branch
from which we bifurcate. A formal first-order calculation in the same
subsection explains the choice of the Dirichlet resonance and identifies the
deformation direction. In \S\ref{sec:conformal} we parametrise
nearby domains by conformal maps and pull back the free-boundary
problem to the unit disk. We show that the resulting equation
$\Fcal(p,f;\kappa)=0$ is real analytic and equivalent to the original
problem. In \S\ref{sec:bifurcation} we compute the linearisation along
the trivial branch. The correction $B_\kappa f$ allows us to reduce
the linearised problem to a Dirichlet problem for the field variation,
with the boundary equation then determining the deformation. At
$\kappa=\kappa_*$, we verify the assumptions of the analytic Crandall--Rabinowitz
theorem to obtain a real-analytic branch of exact non-scattering
configurations, which we parametrise by a dilation-invariant geometric
quantity $\eps$.

In \S\ref{sec:geometry} we prove that the domains on the branch are
strictly convex, have real-analytic boundaries, and are non-circular away
from the starting disk. We normalise the mean polar radius to one and
use rotational symmetry to show that the wave numbers and material coefficient are
even in the parameter $\eps$. In \S\ref{sec:coefficients} we
use exact integral identities to compute the leading corrections to the
polar radius and material coefficient. Comparing the second-order shape correction
with that of a normalised ellipse rules out ellipses and completes the
proof of Theorem~\ref{thm:non-scattering-branch}.

We defer the auxiliary technical tools to the Appendix for clarity of presentation. More precisely, we state the analytic Crandall--Rabinowitz theorem in
\S\ref{A2} and prove the Bessel identities used in the expansions
in Lemma~\ref{lem:bessel-data}.

%%%%%%%%%%%%%%%%%%%%%%%%%%%%%%%%%%%%%%%%%%%%%%%%%%%%%%%%%%%%%%%%
\section{Construction of non-scattering analytic domains}\label{main}
In this section we construct a family of domains $\widetilde\Omega_\eps$
with spatially constant coefficients $\sigma_\eps$. Each domain
$\widetilde\Omega_\eps$ is non-scattering for the radial Herglotz incident wave
$v_{\tilde k_\eps}$ defined by \eqref{eq:incident}, where $\tilde k_\eps$ is the
exterior wave number. The construction proceeds by bifurcation from
the resonant disk $\widetilde\Omega_0=\D$, with parameters
$\tilde k_0=k_*$, $\tilde\kappa_0=\kappa_*$, and $\sigma_0=\sigma_*$, which we describe next.
The family $\Omega_\eps$ of Theorem~\ref{thm:non-scattering-branch} is obtained
from $\widetilde\Omega_\eps$ by a dilation in \S\ref{sec:geometry}.

\subsection{The resonant starting disk}\label{sec:disk}
Note that for any $t>0$ such that $J_0(t)\ne 0$,  $d_0(t)$ defined  by  \eqref{eq:diskmatch} is the eigenvalue of the Dirichlet-to-Neumann map on $\partial\D$ for the Helmholtz equation $(\Delta+t^2)u=0$ in $\D$, corresponding to the constant eigenfunction, see, e.g., \cite[Chapter 7]{LMP}. It is well known that it is monotone decreasing on any interval of the positive real line not containing zeros of $J_0$. Indeed,  for $t>0$ with $J_0(t)\ne0$, Bessel's equation gives
\begin{equation}\label{eq:d0-derivative}
d_0'(t)=-t-\frac{d_0(t)^2}{t}<0.
\end{equation}
The next lemma shows that the parameters in \eqref{eq:seeds} are well defined.
\begin{lemma}\label{lem1}
For $\kappa_*=j_{2,1}$, we have $J_0(\kappa_*)\ne0$ and
\begin{equation}\label{eq:seed-d0}
d_0(\kappa_*)=-\frac{\kappa_*^2}{2}.
\end{equation}
The matching equation $d_0(k_*)=d_0(\kappa_*)$ in \eqref{eq:seeds}
has a unique solution $k_*\in(0,j_{0,1})$.
\end{lemma}

\begin{proof} 
Since $J_2(\kappa_*)=0$, the relation $J_2(t)=2J_1(t)/t-J_0(t)$ gives $J_1(\kappa_*)=\kappa_*J_0(\kappa_*)/2$. Hence $J_0(\kappa_*)\ne0$, as $J_0$ and $J_0'=-J_1$ have no common zeros, and $d_0(t)=-tJ_1(t)/J_0(t)$ proves  \eqref{eq:seed-d0}.
As $d_0(t)$ is strictly monotone decreasing on the interval $\left(0, j_{0,1}\right)$ from $0$ to $-\infty$, it takes every negative value in this interval exactly once, including $d_0(\kappa_*)$. 
\end{proof}

Since $d_0'(k_*)\ne0$, the analytic implicit function theorem gives an
open interval $I\ni\kappa_*$ and an analytic function $k=k(\kappa)$ on
$I$ with $k(\kappa_*)=k_*$ and $d_0(k(\kappa))=d_0(\kappa)$. In what
follows we shrink $I$ whenever necessary without further mention. Throughout \S\ref{main}, we write $k=k(\kappa)$. Define
\begin{equation}\label{eq:w0}
w_0(r;\kappa)=
\frac{J_0(k)}{J_0(\kappa)}J_0(\kappa r)-J_0(kr),
\qquad
S_\kappa(x)=(k^2-\kappa^2)J_0(k|x|).
\end{equation}
By \eqref{eq:disk-fields}--\eqref{eq:diskmatch},
\begin{equation}\label{eq:diskcurve}
(\Delta+\kappa^2)w_0=S_\kappa\quad\text{in }\D,\qquad
w_0=\partial_rw_0=0\quad\text{on }\mathbb S^1.
\end{equation}
This is the \emph{trivial branch}: the domain stays circular, although
the fields are nonzero. Evaluating the radial  equation 
\eqref{eq:diskcurve} at $r=1$, where $w_0$ and $\partial_rw_0$ vanish, we
obtain 
\begin{equation}\label{eq:Ckappa}
C_\kappa:=\partial_r^2w_0(1;\kappa)
=(k^2-\kappa^2)J_0(k)\ne0.
\end{equation}

The choice $\kappa_*=j_{2,1}$ is motivated by the following formal
first-order calculation, which asks whether the disk and its radial field
$w_0$ can be deformed while preserving both zero boundary conditions in
\eqref{eq:clamped}. Fix $\kappa\in I$ and $k=k(\kappa)$, and suppose that for
a small real parameter $\tau$ the problem \eqref{eq:clamped} has a solution
$w_\tau$ in a domain with boundary $r=1+\tau \chi(\theta)+O(\tau^2)$, where
$w_\tau=w_0+\tau p+O(\tau^2)$ at fixed interior points of $\D$; thus $p$ is
the Eulerian variation of the field. Since the wave numbers are fixed, the
source does not depend on $\tau$. Differentiating the equation and the two
boundary conditions at $\tau=0$, and using that $\nabla w_0=0$ and
$\partial_r^2w_0=C_\kappa$ on $\mathbb S^1$, we obtain
\begin{equation}\label{eq:formal-linearisation}
\begin{cases}
(\Delta+\kappa^2)p=0 & \text{in }\D,\\
p=0 & \text{on }\mathbb S^1,\\
\partial_rp+C_\kappa \chi=0 & \text{on }\mathbb S^1.
\end{cases}
\end{equation}
If $\kappa^2$ is not a Dirichlet eigenvalue of the disk, then $p=0$ and,
as $C_\kappa\ne0$, also $\chi=0$: there is no non-trivial first-order
deformation. Displacements $\chi$ of angular order zero and one correspond, to
first order, to a dilation and a translation of the disk, so angular order
two is the lowest that changes its shape. Accordingly, we take
$\kappa_*=j_{2,1}$, for which
\begin{equation}\label{eq:psi}
\psi(r,\theta)=J_2(\kappa_*r)\cos2\theta
\end{equation}
is a Dirichlet eigenfunction. We set
\begin{equation}\label{eq:constants}
C_*=C_{\kappa_*},\qquad
\beta_*=\kappa_*J_2'(\kappa_*)\ne0,
\end{equation}
where $\beta_*\ne0$ because the positive zeros of $J_2$ are simple. The
choice $p=-(C_*/\beta_*)\psi$ solves \eqref{eq:formal-linearisation}
with $\chi(\theta)=\cos2\theta$, so that the boundary begins to move as
$r=1+\tau\cos2\theta+O(\tau^2)$. The companion mode $\sin2\theta$ is a rotation
of this one; we remove it by working with domains invariant under
reflection in both coordinate axes.

%%%%%%%%%%%%%%%%%%%%%%%%%%%%%%%%%%%%%%%%%%%%%%%%%%%%%%%%%%%%%%%%
\subsection{A conformal formulation of non-scattering on the unit disk}\label{sec:conformal}
A key point of our approach is to conformally  transplant \eqref{eq:clamped}   to the fixed disk $\D$. To this end, identifying $\R^2$ with $\C$, we define the real
Banach space
\[
\Hcal
=
\left\{
f\in C^{4,\alpha}(\overline{\D};\C):
f\text{ is holomorphic in }\D,\quad
f(-z)=-f(z),\quad
f(\overline z)=\overline{f(z)}
\right\}.
\]
Thus the elements of $\Hcal$ are the functions
\begin{equation}\label{eq:ftaylor}
f(z)=a_1z+a_3z^3+a_5z^5+\cdots,
\qquad a_{2j+1}\in\R,
\end{equation}
of the stated regularity. Although the solutions we construct are real analytic, the Banach space framework is essential to our proof.

For $f$ of sufficiently small norm in $\Hcal$, we set
\begin{equation}\label{eq:conformal-map}
\Phi_f(z)=z+f(z),
\qquad
\Omega_f=\Phi_f(\D).
\end{equation}
In particular, we take $f$ small enough that
$\|f'\|_{L^\infty(\D)}<1$. For $z_1,z_2\in\overline{\D}$,
\[
|\Phi_f(z_1)-\Phi_f(z_2)|
\ge
\left(1-\|f'\|_{L^\infty(\D)}\right)|z_1-z_2|.
\]
Thus $\Phi_f$ is injective on $\overline{\D}$ and conformal in $\D$,
with $\Phi_f(0)=0$ and $\Phi_f'(0)>0$. The symmetries of $f$ make
$\Omega_f$ invariant under reflection in both coordinate axes.

Let $C_{\mathrm{even},\pi}^{4,\alpha}(\mathbb S^1;\R)$ denote the
subspace of $C^{4,\alpha}(\mathbb S^1;\R)$ consisting of functions
that are even and $\pi$-periodic. Their Fourier series consist of
$\cos(2j\theta)$ modes, $j\ge0$. For
$f(z)=\sum_{j\ge0}a_{2j+1}z^{2j+1}$, the normal component of the
deformation on $\mathbb S^1$ is
\begin{equation}\label{eq:Tmeaning}
Tf(\theta):=\operatorname{Re}\left(e^{-\ir\theta}f(e^{\ir\theta})\right)=\sum_{j\ge0}a_{2j+1}\cos(2j\theta).
\end{equation}
In particular, $T(z^3)=\cos 2\theta$ gives the displacement found
in \S\ref{sec:disk}, and $T(z)=1$ corresponds to an
infinitesimal dilation. We retain the dilation as an unknown until
the scale is normalised in \S\ref{sec:geometry}.

\begin{lemma}\label{lem:boundary-isomorphism}
The operator
$T:\Hcal\to C_{\mathrm{even},\pi}^{4,\alpha}(\mathbb S^1;\R)$
is a bounded real-linear isomorphism with inverse
\begin{equation}\label{eq:Tinv}
T^{-1}\left(\sum_{j\ge0}h_j\cos(2j\theta)\right)(z)
=
z\sum_{j\ge0}h_jz^{2j}.
\end{equation}
\end{lemma}

\begin{proof}
Boundedness and injectivity follow from \eqref{eq:Tmeaning}.
Given $h\in C_{\mathrm{even},\pi}^{4,\alpha}(\mathbb S^1;\R)$, let $U$ be its harmonic extension
to $\D$ and let $V$ be the harmonic conjugate of $U$ normalised by
$V(0)=0$. Set $H=U+\ir V$. By the regularity up to the boundary, the higher-order 
global Schauder estimate for the Dirichlet problem on $\D$ (see
\cite[Theorem~6.19 and Problem 6.2]{GilbargTrudinger}) and the maximum principle, we have
$U\in C^{4,\alpha}(\overline{\D})$ with norm bounded by
$C\|h\|_{C^{4,\alpha}(\mathbb S^1)}$. Since
$H'=\partial_1U-\ir\,\partial_2U$ and $H(0)=U(0)$, it follows that
\[
\|H\|_{C^{4,\alpha}(\overline{\D})}
\le C\|h\|_{C^{4,\alpha}(\mathbb S^1)}.
\]
The symmetries of $h$ and the normalisation of $V$ imply
$H(-z)=H(z)$ and $H(\overline z)=\overline{H(z)}$.
It follows that $f(z)=zH(z)$ belongs to $\Hcal$ and
$Tf=\operatorname{Re}H|_{\mathbb S^1}=h$. Expressing $f$ in terms
of the Fourier coefficients of $h$ gives \eqref{eq:Tinv}.
\end{proof}

Having established the conformal parametrisation of the domain, we now
analyse the pullback and lift of the fields between $\Omega_f$ and $\D$. Fix $\kappa\in I$ and $f\in\Hcal$, and write $k=k(\kappa)$. Let $w=u-v_k$ be the difference between the interior total field $u$
and the incident field $v_k$ on $\Omega_f$, and let
$W=w\circ\Phi_f$ be its pullback to $\D$. For a small real parameter $t$, set
$\Phi_{tf}(z)=z+tf(z)$ and $\Omega_{tf}=\Phi_{tf}(\D)$.
At $t=0$, the map is the identity and the domain is $\D$.
For a differentiable family $w_t$ on $\Omega_{tf}$ starting from
the radial field $w_0(\cdot;\kappa)$, the chain rule gives  
\[
\left.\frac{\dr}{\dr t}(w_t(\Phi_{tf}(z)))\right|_{t=0}
=
\left.\frac{\partial w_t(z)}{\partial t}\right|_{t=0}
+\nabla w_0\cdot f(z).
\]
On the right-hand side, $\partial w_t/\partial t$ is taken at a fixed physical
point. Here and below, we identify $f$ with the vector field
$(\operatorname{Re}f,\operatorname{Im}f)$ in expressions involving
real gradients. The term $\nabla w_0\cdot f$ is the first-order
contribution of the change of coordinates. We split  it off by setting
\begin{equation}\label{eq:lift}
B_\kappa f
:=
\nabla w_0(\,\cdot\,;\kappa)\cdot f,
\qquad
W=w_0+p+B_\kappa f,
\end{equation}
and take $(p,f)$ as the unknowns. To first order, $p$ represents the
Eulerian variation introduced in \S\ref{sec:disk}. We choose $p\in\Pcal$, where
\begin{equation}\label{space}
\Pcal:=
\left\{
p\in C^{5,\alpha}(\overline{\D};\R):
p|_{\mathbb S^1}=0,\quad
p(-z)=p(z),\quad
p(\overline z)=p(z)
\right\}.
\end{equation}
Only the Dirichlet condition is built into $\Pcal$; the Neumann condition
will be part of the nonlinear equation. The lift $B_\kappa f$ gains one
derivative, which puts it in the same H\"older space as $p$; this is made
precise in the lemma below.

\begin{lemma}\label{lem:lift-regularity}
Let $I$ be the interval introduced in \S\ref{sec:disk}. For every $\kappa\in I$, the operator
$B_\kappa:\Hcal\to\Pcal$ is bounded, with
\[
\|B_\kappa f\|_{C^{5,\alpha}(\overline{\D})}
\le C\|f\|_{C^{4,\alpha}(\overline{\D})}
\]
for a constant $C$ independent of $\kappa\in I$. Moreover,
$\kappa\mapsto B_\kappa$ is real analytic from $I$ to
$\mathcal L(\Hcal,\Pcal)$.
\end{lemma}

\begin{proof}
By \eqref{eq:diskcurve}, $\nabla w_0=0$ on
$\mathbb S^1$, since $w_0$ is radial and $\partial_rw_0(1;\kappa)=0$. Hence $B_\kappa f$ has zero trace, and its reflection
symmetries follow from those of $w_0$ and $f$. Write $f=(f_1,f_2)$.
The components of $f$ are harmonic, so
\begin{equation}\label{eq:liftlap}
\Delta(B_\kappa f)
=
\nabla(\Delta w_0)\cdot f
+
2\sum_{i,j=1}^2
(\partial_{ij}w_0)(\partial_i f_j).
\end{equation}
After shrinking $I$ if necessary, the right-hand side has
$C^{3,\alpha}(\overline{\D})$ norm bounded by
$C\|f\|_{C^{4,\alpha}(\overline{\D})}$ uniformly in $\kappa\in I$.
The regularity up to the boundary and the higher-order global Schauder estimate for the
Dirichlet problem on $\D$ \cite[Theorem~6.19 and Problem 6.2]{GilbargTrudinger} give
the stated bound. By \eqref{eq:w0}, the coefficients in \eqref{eq:liftlap} depend real
analytically on $\kappa$ in the corresponding H\"older norms.
Applying the bounded Dirichlet inverse of the Laplacian gives the
real-analytic dependence of $B_\kappa$ in operator norm.
\end{proof}

We are now ready to use the conformal map $\Phi_f$ to reformulate
the non-scattering problem \eqref{eq:clamped} on $\Omega_f$ as an equivalent nonlinear
problem on the unit disk.

Set 
\begin{equation}\label{space2}
\begin{split}
\Xcal&:=\Pcal\times\Hcal,\\
\Ycal&:=\left\{G\in C^{3,\alpha}(\overline{\D};\R):
G(-z)=G(z),\ G(\overline{z})=G(z)\right\}
\times C_{\mathrm{even},\pi}^{4,\alpha}(\mathbb S^1;\R).
\end{split}
\end{equation}
The first factor of $\Ycal$ has unrestricted boundary trace. For $f$ sufficiently small, and with $W$ given by \eqref{eq:lift}, define
\begin{equation}\label{eq:F}
\Fcal(p,f;\kappa)
=
\begin{pmatrix}
\Delta W+\kappa^2|\Phi_f'|^2W-|\Phi_f'|^2(S_\kappa\circ\Phi_f)\\[2mm]
\partial_r W|_{\mathbb S^1}
\end{pmatrix},
\end{equation}
where $S_\kappa$ is given by \eqref{eq:w0}.

\begin{proposition}\label{prop:fixed-disk-formulation}
There exists a neighbourhood $\mathcal U$ of $(0,0)$ in $\Xcal$
such that $\Fcal:\mathcal U\times I\to\Ycal$ is well defined
and real analytic, and
\begin{equation}\label{eq:trivial}
\Fcal(0,0;\kappa)=0
\qquad\text{for every }\kappa\in I.
\end{equation}
For $(p,f;\kappa)\in\mathcal U\times I$, the equation
$\Fcal(p,f;\kappa)=0$ is equivalent to
$w=(w_0+p+B_\kappa f)\circ\Phi_f^{-1}$ satisfying
\begin{equation}\label{eq:clamped-f}
\begin{cases}
(\Delta+\kappa^2)w=S_\kappa
&\text{in }\Omega_f,\\
w=\partial_\nu w=0
&\text{on }\partial\Omega_f,
\end{cases}
\end{equation}
where $\nu$ is the outward unit normal to $\partial\Omega_f$.
\end{proposition}

\begin{proof}
Choose $\mathcal U$ small enough that $\Phi_f$ is conformal for every
$(p,f)\in\mathcal U$. In the decomposition \eqref{eq:lift}, $w_0$ is
smooth on $\overline{\D}$, $p\in\Pcal$ by definition, and
$B_\kappa f\in\Pcal$ by Lemma~\ref{lem:lift-regularity}. It follows
that $W\in C^{5,\alpha}(\overline{\D})$, so that
\[
\Delta W\in C^{3,\alpha}(\overline{\D}),
\qquad
\partial_r W|_{\mathbb S^1}\in C^{4,\alpha}(\mathbb S^1).
\]
Since $|\Phi_f'|^2\in C^{3,\alpha}(\overline{\D})$ and both components in
\eqref{eq:F} have the prescribed symmetries, $\Fcal$ takes values
in $\Ycal$. The map $f\mapsto|\Phi_f'|^2=|1+f'|^2$ is polynomial in the real and imaginary
parts of $f'$. By \eqref{eq:w0} and Lemma~\ref{lem:lift-regularity},
$W$ depends real analytically on $(p,f;\kappa)$. The entire power
series for $J_0(k|y|)$ in $k^2(y_1^2+y_2^2)$ gives the real-analytic
dependence of $S_\kappa\circ\Phi_f$ on $(f,\kappa)$. Continuity of
multiplication in H\"older spaces now gives the real analyticity
of $\Fcal$. At $(p,f)=(0,0)$, we have $W=w_0$,
$\Phi_0=\id$ and $|\Phi_0'|^2=1$, so \eqref{eq:trivial} follows
from \eqref{eq:diskcurve}.

Set $w=W\circ\Phi_f^{-1}$. By conformal covariance of the Laplacian,
\[
\Delta W=|\Phi_f'|^2(\Delta w)\circ\Phi_f,
\]
therefore the first component of $\Fcal$ is
\[
|\Phi_f'|^2\left((\Delta+\kappa^2)w-S_\kappa\right)\circ\Phi_f.
\]
Since $|\Phi_f'|^2>0$, this component vanishes precisely when the interior
equation in \eqref{eq:clamped-f} holds. Each term in the decomposition
of $W$ has zero trace on $\mathbb S^1$, so $w=0$ on
$\partial\Omega_f$.

At $e^{\ir\theta}\in\mathbb S^1$, the differential of $\Phi_f$
maps the outward radial unit vector $e_r$ to
$|\Phi_f'(e^{\ir\theta})|\nu$ at the corresponding point of
$\partial\Omega_f$. Hence
\[
\partial_r W(e^{\ir\theta})
=
\left|\Phi_f'(e^{\ir\theta})\right|\,
(\partial_\nu w)(\Phi_f(e^{\ir\theta})).
\]
The factor $|\Phi_f'(e^{\ir\theta})|$ is positive, so the second
component of $\Fcal$ vanishes precisely when $\partial_\nu w=0$
on $\partial\Omega_f$.
\end{proof}

%%%%%%%%%%%%%%%%%%%%%%%%%%%%%%%%%%%%%%%%%%%%%%%%%%%%%%%%%%%%%%%%
\subsection{Bifurcation from the resonant disk}
\label{sec:bifurcation}
We now prove the existence assertion in
Theorem~\ref{thm:non-scattering-branch} by constructing a branch of
non-scattering domains that bifurcates from the resonant disk constructed in \S\ref{sec:disk}.
The basic tool is the analytic form of the Crandall--Rabinowitz theorem
\cite[Theorem~8.3.1]{BT}, stated in \S\ref{A2} (see \cite[Theorem~1.7]{CR}
for the classical formulation), which we apply to the equation
$\Fcal(p,f;\kappa)=0$ near $(0,0;\kappa_*)$.

Proposition~\ref{prop:fixed-disk-formulation} gives the real
analyticity of $\Fcal$ and the trivial branch
$\Fcal(0,0;\kappa)=0$ needed to apply this theorem. The remaining assumptions concern
the linearisation along this branch at $\kappa=\kappa_*$. The next lemma gives the Fr\'echet derivative of $\Fcal(p,f;\kappa)$ with respect to $(p,f)$ at $(0,0;\kappa)$.
Recall that $\Xcal$ and $\Ycal$ are given by \eqref{space2}.
\begin{lemma}\label{prop:linearization}
For every $\kappa\in I$, the operator
$L_\kappa:=D_{(p,f)}\Fcal(0,0;\kappa):\Xcal\to\Ycal$
is given by
\begin{equation}\label{eq:L}
L_\kappa(q,g)
=
\begin{pmatrix}
(\Delta+\kappa^2)q\\[1mm]
\left.\partial_rq\right|_{\mathbb S^1}+C_\kappa Tg
\end{pmatrix},
\end{equation}
where $C_\kappa$ and $T$ are defined in \eqref{eq:Ckappa} and
\eqref{eq:Tmeaning}, respectively.
\end{lemma}

\begin{proof}
Fix $\kappa\in I$ and $(q,g)\in\Xcal$. For small real $t$, take
$(p,f)=t(q,g)$, so that
$W=w_0+t(q+B_\kappa g)$ with $W$ given by \eqref{eq:lift}. Since $|\Phi_{tg}'|^2=|1+tg'|^2$, we have
\[
\left.\frac{\dr}{\dr t}|\Phi_{tg}'|^2\right|_{t=0}
=2\operatorname{Re}g',
\qquad
\left.\frac{\dr}{\dr t}(S_\kappa\circ\Phi_{tg})\right|_{t=0}
=\nabla S_\kappa\cdot g.
\]
Hence, differentiating the first component of $\Fcal$ gives
\begin{equation}\label{eq:linearization-first}
\left[L_\kappa(q,g)\right]_1
=
(\Delta+\kappa^2)(q+B_\kappa g)
+2\operatorname{Re}g'\left(\kappa^2w_0-S_\kappa\right)
-\nabla S_\kappa\cdot g,
\end{equation}
where the subscript indicates the component. Let $D$ and $D^2$ denote the Jacobian and Hessian, respectively, and identify $\C$ with $\R^2$. By the Cauchy--Riemann equations, the symmetric part of $Dg$, that is,
$(Dg+(Dg)^\top)/2$, is $(\operatorname{Re}g')\,\id$. Since $D^2w_0$ is symmetric, we have
\[
\sum_{i,j=1}^2(\partial_{ij}w_0)(\partial_i g_j)
=(\operatorname{Re}g')\Delta w_0.
\]
Thus \eqref{eq:liftlap} and \eqref{eq:diskcurve} give
\begin{equation}\label{eq:lift-helmholtz}
(\Delta+\kappa^2)B_\kappa g = \nabla\left((\Delta+\kappa^2)w_0\right)\cdot g
+2\operatorname{Re}g'\,\Delta w_0 = \nabla S_\kappa\cdot g +2\operatorname{Re}g'\,\Delta w_0.
\end{equation}
Substituting this identity into \eqref{eq:linearization-first},
we obtain
\[
\left[L_\kappa(q,g)\right]_1
=
(\Delta+\kappa^2)q
+2\operatorname{Re}g'
  \left((\Delta+\kappa^2)w_0-S_\kappa\right)
=
(\Delta+\kappa^2)q.
\]
For the boundary component, the fact that $w_0$ is a radial function gives
\[
B_\kappa g(re^{\ir\theta})
=
\partial_rw_0(r;\kappa)
\operatorname{Re}\left(e^{-\ir\theta}g(re^{\ir\theta})\right).
\]
Differentiating in $r$ at $r=1$ and using
$\partial_rw_0(1;\kappa)=0$, we find
$\partial_r(B_\kappa g)|_{\mathbb S^1}=C_\kappa Tg$.
This gives the second component in \eqref{eq:L}.
\end{proof}

The contribution of $B_\kappa g$ cancels the other terms involving
$g$ in the interior equation. Since $C_\kappa\ne0$ and $T$ is an
isomorphism by Lemma~\ref{lem:boundary-isomorphism}, we solve
$L_\kappa(q,g)=(G,h)$ by first solving the Dirichlet problem for $q$
and then recovering $g$ from the boundary equation:
\begin{equation}\label{eq:linear-dirichlet}
\begin{aligned}
(\Delta+\kappa^2)q&=G\quad\text{in }\D,
&\qquad q&=0\quad\text{on }\mathbb S^1,\\
g&=C_\kappa^{-1}T^{-1}
\left(h-\left.\partial_rq\right|_{\mathbb S^1}\right).
\end{aligned}
\end{equation}
Thus all solvability conditions come from the Dirichlet problem.
In particular, $L_\kappa$ is an isomorphism whenever $\kappa^2$
lies outside the Dirichlet spectrum in the prescribed
reflection-symmetric class.

Having obtained an explicit expression for the linearised operator, we are now able to verify the remaining hypotheses of the Crandall--Rabinowitz theorem.
Recall that $\psi$, $C_*$ and $\beta_*$ are defined in \eqref{eq:psi} and \eqref{eq:constants}.
\begin{proposition}\label{prop:kernel-range}
The operator $L_{\kappa_*}:\Xcal\to\Ycal$ has the following
properties.
\begin{enumerate}[\rm(i)]
\item Its kernel is one-dimensional,
\begin{equation}\label{eq:kernel}
\ker L_{\kappa_*}
=
\operatorname{span}\{\xi_*\},
\qquad
\xi_*=
\left(-\frac{C_*}{\beta_*}\psi,\ z^3\right).
\end{equation}

\item Its range is closed and has codimension one,
\begin{equation}\label{eq:range}
\Ran L_{\kappa_*}
=
\left\{(G,h)\in\Ycal:
\int_{\D}G\psi\,\dr x=0\right\}.
\end{equation}

\item The transversality condition holds,
\begin{equation}\label{eq:transversality}
(\left.\partial_\kappa L_\kappa\right|_{\kappa=\kappa_*})\xi_*
\notin\Ran L_{\kappa_*}.
\end{equation}
\end{enumerate}
\end{proposition}

\begin{proof}
(i) Let $L_{\kappa_*}(q,g)=0$. In the first component, $q$ solves the homogeneous Dirichlet problem in $\D$.  At $\kappa_*=j_{2,1}$ the Dirichlet eigenspace is spanned by
the two modes $J_2(\kappa_*r)\cos2\theta$ and
$J_2(\kappa_*r)\sin2\theta$. Indeed, $J_0(\kappa_*)\ne0$ by
Lemma~\ref{lem1}, and the recurrence relation for the Bessel functions gives
$J_1(\kappa_*)=(\kappa_*/2)J_0(\kappa_*)\ne0$. By Sturm comparison for the zeros of Bessel's equation, the first positive zero $j_{n,1}$ of $J_n$
is strictly increasing in $n$. Thus $j_{n,1}>j_{2,1}=\kappa_*$
for $n\ge3$, so $J_n(\kappa_*)\ne0$. The reflection symmetries in $\Pcal$ select the cosine mode,
so $q=c\psi$ for some $c\in\R$. Since
$\partial_r\psi(1,\theta)=\beta_*\cos2\theta$ and
$T(z^3)=\cos2\theta$, \eqref{eq:linear-dirichlet} with
$(G,h)=(0,0)$ gives $g=-c(\beta_*/C_*)z^3$. This proves
\eqref{eq:kernel}.

(ii) By the Fredholm alternative, the Dirichlet problem in
\eqref{eq:linear-dirichlet} at $\kappa=\kappa_*$ is solvable
if and only if $G$ is orthogonal to both Dirichlet modes. Since
$G$ is even in $\theta$, its inner product with the sine mode
vanishes. The remaining condition is
\[
\ell(G,h):=\int_{\D}G\psi\,\dr x=0.
\]

When this condition holds, we average a solution over the group
generated by the two coordinate reflections. Since $G$ has these
symmetries, the average solves the same equation and, by Schauder
regularity, belongs to $\Pcal$. For this $q$, the boundary equation
determines $g$ uniquely for each $h$, since $C_*\ne0$ and $T$ is an
isomorphism. Thus $\Ran L_{\kappa_*}=\ker\ell$, which proves
\eqref{eq:range}. The functional $\ell$ is continuous and nonzero,
since $\ell(\psi,0)=\int_{\D}\psi^2\,\dr x>0$.

(iii) By \eqref{eq:range}, the range condition $\ell(G,h)=0$ depends only on $G$.
By \eqref{eq:L}, the $\kappa$-derivative of the first component of
$L_\kappa(q,g)$, with $(q,g)$ fixed, is $2\kappa q$. Taking
$(q,g)=\xi_*$ as in \eqref{eq:kernel}, we obtain
\[
\ell\left(
\left.\partial_\kappa L_\kappa\right|_{\kappa=\kappa_*}\xi_*
\right)
=
-2\kappa_*\frac{C_*}{\beta_*}\int_{\D}\psi^2\,\dr x
\ne0.
\]
Since $\Ran L_{\kappa_*}=\ker\ell$ by \eqref{eq:range}, this
proves \eqref{eq:transversality}.
\end{proof}

Note that the kernel of $L_{\kappa_*}$ is one-dimensional because the
reflection symmetries select the cosine mode from the two-dimensional
Dirichlet eigenspace. Since the range is closed and has codimension one
by \eqref{eq:range}, $L_{\kappa_*}$ is Fredholm of index zero.

We now apply the analytic bifurcation theorem of Crandall and
Rabinowitz to construct the family of non-scattering domains in
Theorem~\ref{thm:non-scattering-branch}. We parametrise the branch by a
geometric quantity which is invariant under dilations: the
$\cos2\theta$ Fourier coefficient of the polar radius of the domain,
divided by the mean polar radius.

\begin{proposition}\label{prop:bifurcating-branch}
There exist $\eps_0>0$ and a real-analytic curve
\[
(-\eps_0,\eps_0)\ni\eps\longmapsto
(\tilde p_\eps,\tilde f_\eps;\tilde\kappa_\eps)\in\Xcal\times I
\]
with $(\tilde p_0,\tilde f_0;\tilde\kappa_0)=(0,0;\kappa_*)$ such that
\begin{equation}\label{eq:branch-tangent}
\Fcal(\tilde p_\eps,\tilde f_\eps;\tilde\kappa_\eps)=0,
\qquad
(\tilde p_\eps,\tilde f_\eps)=\eps\xi_*+O_{\Xcal}(\eps^2).
\end{equation}
Every zero of $\Fcal$ sufficiently close to $(0,0;\kappa_*)$
belongs either to the trivial branch $(p,f)=(0,0)$ or to this curve. Furthermore, set
\begin{equation}\label{eq:phibranch}
\widetilde\Phi_\eps(z)=z+\tilde f_\eps(z)=z+\eps z^3+O_{C^{4,\alpha}}(\eps^2),
\qquad
\widetilde\Omega_\eps=\widetilde\Phi_\eps(\D),
\qquad
\tilde k_\eps=k(\tilde\kappa_\eps).
\end{equation}
Then the domain $\widetilde\Omega_\eps$, with constant coefficient
$\sigma_\eps=\tilde\kappa_\eps^2/\tilde k_\eps^2>1$, is non-scattering for the incident
field $v_{\tilde k_\eps}$. It is star-shaped with respect to the origin, and
its polar radius $\tilde r_\eps$ and mean polar radius
$R_\eps:=\frac1{2\pi}\int_0^{2\pi}\tilde r_\eps(\theta)\,\dr\theta$ satisfy
\begin{equation}\label{fin}
\begin{gathered}
\tilde r_\eps(\theta)=1+\eps\cos 2\theta+O(\eps^2)
\quad\text{uniformly in }\theta,\\
R_\eps=1+O(\eps^2),\qquad
\frac1{\pi R_\eps}\int_0^{2\pi}\tilde r_\eps(\theta)\cos 2\theta\,\dr\theta=\eps.
\end{gathered}
\end{equation}
Moreover, $R_\eps$ is real analytic in $\eps$.
\end{proposition}

\begin{proof}
By Proposition~\ref{prop:fixed-disk-formulation},
Lemma~\ref{prop:linearization} and
Proposition~\ref{prop:kernel-range}, all the hypotheses of
Theorem~\ref{thm:analytic-cr} are satisfied with $F=\Fcal$, $x=(p,f)$,
$\lambda=\kappa$, $\lambda_*=\kappa_*$ and $\phi=\xi_*$. This gives
$s_0>0$ and a real-analytic curve
$(-s_0,s_0)\ni s\mapsto(p_s,f_s;\kappa_s)\in\Xcal\times I$ with
$(p_0,f_0;\kappa_0)=(0,0;\kappa_*)$ such that
\begin{equation}\label{eq:branch-s}
\Fcal(p_s,f_s;\kappa_s)=0,
\qquad
(p_s,f_s)=s\xi_*+O_{\Xcal}(s^2),
\end{equation}
together with the stated local description of the zeros. It remains to
replace $s$ by the parameter $\eps$.

Let $\gamma_s(t)=e^{\ir t}+f_s(e^{\ir t})$, $t\in[0,2\pi]$, which parametrises
the boundary of the domain $\Omega^{(s)}:=(\id+f_s)(\D)$. Since the second component of $\xi_*$ is
$z^3$, \eqref{eq:branch-s} gives
$\gamma_s(t)=e^{\ir t}(1+se^{2\ir t})+O(s^2)$ in $C^{4,\alpha}(\mathbb S^1)$. Write
$\gamma_s(t)=|\gamma_s(t)|e^{\ir\vartheta_s(t)}$, where
$\vartheta_s$ is the continuous argument with $\vartheta_s(0)=0$.
For small $s$, we have
$\vartheta_s'=\operatorname{Im}\frac{\gamma_s'}{\gamma_s}>0$, so the change of variable $\theta=\vartheta_s(t)$
is strictly increasing. Hence $\Omega^{(s)}$ is star-shaped with respect to
the origin, and its boundary has a polar representation $r=r_s(\theta)$.
The change of variables
expresses the mean polar radius $R_s$ of $\Omega^{(s)}$ and the
$\cos2\theta$ Fourier coefficient $E(s)$ of $r_s/R_s$ in terms
of $\gamma_s$ and $\gamma_s'$ as follows:
\begin{align*}
R_s
&:=
\frac1{2\pi}\int_0^{2\pi}r_s(\theta)\,\dr\theta
=
\frac1{2\pi}\int_0^{2\pi}
|\gamma_s|\,
\operatorname{Im}\frac{\gamma_s'}{\gamma_s}\,\dr t,
\\
E(s)
&:=
\frac1{\pi R_s}\int_0^{2\pi}r_s(\theta)\cos 2\theta\,\dr\theta
=
\frac1{\pi R_s}\int_0^{2\pi}
|\gamma_s|\,
\operatorname{Re}\frac{\gamma_s^2}{|\gamma_s|^2}\,
\operatorname{Im}\frac{\gamma_s'}{\gamma_s}\,\dr t.
\end{align*}
The integrands are real-analytic expressions in $\gamma_s$ and
$\gamma_s'$. Since $\gamma_s$ depends real analytically on $s$
and $|\gamma_s|$ is uniformly bounded away from zero, both $R_s$
and $E(s)$ are real analytic in $s$. Moreover,
$|\gamma_s(t)|=1+s\cos 2t+O(s^2)$ and $\vartheta_s(t)=t+O(s)$,
which give
\[
r_s(\theta)=1+s\cos 2\theta+O(s^2)
\quad\text{uniformly in }\theta,
\qquad
R_s=1+O(s^2),
\qquad
E(s)=s+O(s^2).
\]
Since $E(0)=0$ and $E'(0)=1$, we may take $\eps=E(s)$ as a new
parameter, with real-analytic inverse $s=s(\eps)=\eps+O(\eps^2)$ defined
for $|\eps|<\eps_0$. We set
$(\tilde p_\eps,\tilde f_\eps;\tilde\kappa_\eps)
:=(p_{s(\eps)},f_{s(\eps)};\kappa_{s(\eps)})$, so that
$\widetilde\Omega_\eps=\Omega^{(s(\eps))}$, $\tilde r_\eps=r_{s(\eps)}$ and
$R_\eps=R_{s(\eps)}$. Then \eqref{eq:branch-tangent}, the expansion
in \eqref{eq:phibranch} and \eqref{fin} follow from the corresponding
statements in the parameter $s$, and the local description of the zeros is
unchanged.

By Proposition~\ref{prop:fixed-disk-formulation}, the field
\[
\tilde w_\eps=
\left(w_0(\,\cdot\,;\tilde\kappa_\eps)+\tilde p_\eps+B_{\tilde\kappa_\eps}\tilde f_\eps\right)
\circ\widetilde\Phi_\eps^{-1}
\]
satisfies \eqref{eq:clamped-f} on $\widetilde\Omega_\eps$ with $\kappa=\tilde\kappa_\eps$.
Thus, setting $\tilde u_\eps=v_{\tilde k_\eps}+\tilde w_\eps$, we find that $(\tilde u_\eps,v_{\tilde k_\eps})$
satisfies \eqref{eq:transmission} on $\widetilde\Omega_\eps$ with wave numbers
$\tilde k_\eps$ and $\tilde\kappa_\eps$. Since $\sigma_\eps$ converges to $\sigma_*>1$ as $\eps\to 0$, we have
$\sigma_\eps>1$ after reducing $\eps_0$ if necessary.
\end{proof}

The local uniqueness in Proposition~\ref{prop:bifurcating-branch}
is within the prescribed reflection-symmetric class with
$k=k(\kappa)$. The next section establishes the geometric properties
of $\widetilde\Omega_\eps$, including real analyticity of the boundary, and
normalises the scale to obtain the family in
Theorem~\ref{thm:non-scattering-branch}.

%%%%%%%%%%%%%%%%%%%%%%%%%%%%%%%%%%%%%%%%%%%%%%%%%%%%%%%%%%%%%%%%%
\section{The non-scattering family: geometry and expansions}
\label{sec:family}

\subsection{Geometry of the branch and its normalisation}
\label{sec:geometry}

Proposition~\ref{prop:bifurcating-branch} gives a branch of domains
$\widetilde\Omega_\eps=\widetilde\Phi_\eps(\D)$ bifurcating from the resonant disk. Each domain
is non-scattering for the incident wave $v_{\tilde k_\eps}$ with the constant
coefficient $\sigma_\eps=\tilde\kappa_\eps^2/\tilde k_\eps^2$. We now establish the geometric
properties of these domains and use dilation invariance to normalise
their mean polar radius to one. Since the parameter $\eps$ is invariant
under dilations, it remains the $\cos 2\theta$ coefficient of the
normalised polar radius, as in Proposition~\ref{prop:normalised-expansions}.

\begin{proposition}\label{prop:branch-geometry}
For all sufficiently small $\eps$, the domain $\widetilde\Omega_\eps$ is bounded,
simply connected, contains the origin, and is invariant under
reflection in both coordinate axes. Its boundary is real analytic
and has everywhere positive curvature. For $\eps\ne0$, the domain
$\widetilde\Omega_\eps$ is non-circular.
\end{proposition}

\begin{proof}
For small $\eps$, we have $\|\tilde f_\eps'\|_{L^\infty(\D)}<1$. Thus the assertions in the first
sentence follow from the properties of $\Phi_f$ established in
\S\ref{sec:conformal}.

Let $\tilde\gamma_\eps(t)=\widetilde\Phi_\eps(e^{\ir t})$. The signed curvature of this
anticlockwise parametrisation is
\[
\widetilde{\mathcal K}_\eps
=
\frac{\operatorname{Im}
\left(\overline{\tilde\gamma_\eps'}\,\tilde\gamma_\eps''\right)}
{|\tilde\gamma_\eps'|^3}.
\]
Since $\tilde\gamma_\eps\to e^{\ir t}$ in $C^{4,\alpha}(\mathbb S^1)$, hence in $C^2(\mathbb S^1)$, as $\eps\to0$, the curvature converges
uniformly to one. Hence $\widetilde{\mathcal K}_\eps>0$ for small $\eps$, and
$\widetilde\Omega_\eps$ is strictly convex.

If $\widetilde\Omega_\eps$ were a disk, its reflection symmetries would force
its centre to be the origin, and its polar radius would be constant.
By the last identity in \eqref{fin}, this is impossible for $\eps\ne0$.
Thus $\widetilde\Omega_\eps$ is non-circular for such $\eps$.

Finally, as $\eps\to 0$ we have $\sigma_\eps\to\sigma_*>1$,  $\tilde k_\eps\to k_*<j_{0,1}$, and
$\widetilde\Phi_\eps\to\id$ uniformly on $\overline{\D}$. Consequently,
$\sigma_\eps\ne1$ and $\tilde k_\eps|x|<j_{0,1}$ on $\partial\widetilde\Omega_\eps$ for
small $\eps$, so $J_0(\tilde k_\eps|x|)$ is positive there. Since $\widetilde\Omega_\eps$
is non-scattering, the regularity theorem in
\cite[Theorem~2.1]{CV} applies at every boundary point and gives
the real analyticity of $\partial\widetilde\Omega_\eps$.
\end{proof}

For the \emph{normalisation} of the branch, we use the dilation invariance of
\eqref{eq:transmission}. If $(u,v_k)$ solves the problem on $\Omega$
with wave numbers $(k,\kappa)$, then for any $R>0$ the pair
$(u(R\,\cdot\,),v_{Rk})$ solves it on $R^{-1}\Omega$ with wave
numbers $(Rk,R\kappa)$. The coefficient $\sigma=\kappa^2/k^2$ is unchanged. We
use this to fix the scale, which was not normalised in the bifurcation
argument. Taking $R=R_\eps$, the mean polar radius of
$\widetilde\Omega_\eps$ from Proposition~\ref{prop:bifurcating-branch}, we define
\begin{equation}\label{eq:normalized-family}
\Phi_\eps=R_\eps^{-1}\widetilde\Phi_\eps,
\qquad
\Omega_\eps=\Phi_\eps(\D)=R_\eps^{-1}\widetilde\Omega_\eps,
\qquad
k_\eps=R_\eps\tilde k_\eps,
\qquad
\kappa_\eps=R_\eps\tilde\kappa_\eps.
\end{equation}
Since $R_\eps$ is real analytic in $\eps$, the maps $\Phi_\eps$ and the wave numbers depend real analytically
on $\eps$. By dilation invariance, $\Omega_\eps$ is non-scattering
for $v_{k_\eps}$ with the same coefficient $\sigma_\eps=\kappa_\eps^2/k_\eps^2$,
and retains the properties in Proposition~\ref{prop:branch-geometry}.
By \eqref{fin}, its polar radius $\rho_\eps=\tilde r_\eps/R_\eps$ satisfies
\eqref{eq:polar-normalisation}. Note that $k_\eps$ and $\kappa_\eps$ need no longer satisfy
$d_0(k_\eps)=d_0(\kappa_\eps)$: that relation only fixed the reference
curve for the bifurcation argument.

The \emph{evenness} of the wave numbers in $\eps$ follows from a rotation argument. Rotation through $\pi/2$ preserves the
equations, the symmetry classes
$\Xcal$ and $\Ycal$, and the relation $k=k(\kappa)$. For $\eps\ne 0$, it maps the zero
$(\tilde p_\eps,\tilde f_\eps;\tilde\kappa_\eps)$ of $\Fcal$ to another non-trivial zero with the same
$\kappa$, which by the uniqueness statement in
Proposition~\ref{prop:bifurcating-branch} lies on the same branch. The
rotation preserves the mean polar radius and changes the sign of the
$\cos2\theta$ coefficient of the polar radius, so by the last identity in
\eqref{fin} the rotated zero corresponds to the parameter $-\eps$. Thus
$\tilde\kappa_{-\eps}=\tilde\kappa_\eps$ and $R_{-\eps}=R_\eps$, and consequently
\begin{equation}\label{eveness}
\rho_{-\eps}(\theta)
=
\rho_\eps(\theta+\pi/2),
\qquad
k_{-\eps}=k_\eps,
\qquad
\kappa_{-\eps}=\kappa_\eps,
\end{equation}
so $k_\eps$, $\kappa_\eps$ and $\sigma_\eps$ are even analytic functions
of $\eps$.

This discussion, together with Proposition~\ref{prop:bifurcating-branch},
provides the basis for the proof of Theorem~\ref{thm:non-scattering-branch}
and the first assertion of Proposition~\ref{prop:normalised-expansions}.
We complete the proofs of both results at the end of \S\ref{sec:coefficients}.

%%%%%%%%%%%%%%%%%%%%%%%%%%%%%%%%%%%%%%%%%%%%%%%%%%%%%%%%%%%%%%%%
\subsection{Shape and material coefficient $\eps$-expansions}
\label{sec:coefficients}
It remains to compute the first nonlinear corrections of the normalised
family \eqref{eq:normalized-family}, that is, to establish \eqref{eq:mainshape} and \eqref{eq:main-contrast}. We do this without expanding the transmitted field, using a well-known characterisation of non-scattering by integral identities which involve only the domain and the two wave numbers; we prove it
for the reader's convenience.
\begin{lemma}\label{prop:moment-characterization}
Let $\Omega\subset\R^2$ be a bounded $C^2$ domain
with connected exterior, and let $k,\kappa>0$ with
$k\ne\kappa$. There exists $w\in H_0^2(\Omega)$ satisfying
$(\Delta+\kappa^2)w=(k^2-\kappa^2)J_0(k|x|)$ in $\Omega$
if and only if
\begin{equation}\label{eq:exact-moments}
\int_\Omega J_0(k|x|)e^{-\ir\kappa\omega\cdot x}\,\dr x=0
\qquad\text{for every }\omega\in\mathbb S^1.
\end{equation}
\end{lemma}

This lemma states  that the inhomogeneity $\Omega$ with coefficient $\sigma=\kappa^2/k^2$ does not scatter $J_0(k|x|)$ at the wave number $k$  if and only if \eqref{eq:exact-moments} is satisfied. Note the resemblance between \eqref{eq:exact-moments} and the criterion for failure of the Pompeiu property \cite{HV3}. We remark that this characterisation also remains valid when $\kappa^2$ is a Dirichlet eigenvalue, which is the case for the starting disk.

\begin{proof}[Proof of Lemma~\ref{prop:moment-characterization}]
If such a $w$ exists, Green's identity with the solution
$e^{-\ir\kappa\omega\cdot x}$ of the Helmholtz equation with wave number
$\kappa$, and the two zero traces of $w$, give \eqref{eq:exact-moments}
multiplied by $k^2-\kappa^2\ne0$.

Conversely, suppose that \eqref{eq:exact-moments} holds, and  let
$\frac{\ir}{4}H_0^{(1)}(\kappa|x|)$ be the outgoing fundamental solution of  the Helmholtz equation in $\R^2$ with wave number $\kappa$. Define the volume potential 
\[
w(x)=\frac{\ir}{4}(\kappa^2-k^2)\int_{\Omega} 
H_0^{(1)}(\kappa|x-y|)J_0(k|y|)\,\dr y.
\]
Then $w\in H_{\mathrm{loc}}^2(\R^2)$ and satisfies
$(\Delta+\kappa^2)w=(k^2-\kappa^2)\mathbf1_\Omega J_0(k|x|)$ in $\R^2$. Up to a constant factor, the far-field pattern $w^\infty(\hat x)$ equals the left-hand side of \eqref{eq:exact-moments} with $\omega=\hat x$, hence it is zero.  
Rellich's lemma implies $w=0$ outside $\overline\Omega$ \cite{ColtonKress2019}. Hence both traces of $w$ on
$\partial\Omega$ vanish, and $w|_\Omega\in H_0^2(\Omega)$.
\end{proof}

We derive the shape and material coefficient expansions from the integral identities
\eqref{eq:exact-moments}. We set
\begin{equation}\label{eq:radial-integrals}
b_n(r;k,\kappa):=rJ_0(kr)J_n(\kappa r),
\qquad
F_n(\rho;k,\kappa):=\int_0^\rho b_n(r;k,\kappa)\,\dr r,
\end{equation}
and write $\langle L\rangle=(2\pi)^{-1}\int_0^{2\pi}L(\theta)\,\dr\theta$
for the angular mean. For a domain
$\Omega=\{r(\cos\theta,\sin\theta):0\le r<\rho(\theta)\}$ invariant under
reflection in both coordinate axes, the Jacobi--Anger expansion for the exponential
function expresses \eqref{eq:exact-moments} in terms of angular moments.
The sine moments and the moments of odd order vanish by symmetry, so
\eqref{eq:exact-moments} is equivalent to
\begin{equation}\label{eq:radial-moments}
M_n=M_n(k,\kappa,\rho)
:=
\left\langle\cos(n\theta)F_n(\rho(\theta);k,\kappa)\right\rangle
=0,
\qquad n=0,2,4,\ldots.
\end{equation}
These identities hold along the normalised branch
\eqref{eq:normalized-family}.

Since $\eps\mapsto\Phi_\eps$ is real analytic with values in
$C^{4,\alpha}(\overline{\D};\C)$, the boundary parametrisation of
$\Omega_\eps$ is jointly $C^4$ in the parameter and the boundary
variable. The derivative of the polar angle is bounded away from
zero for small $\eps$, so the implicit function theorem gives joint
$C^4$ regularity of the mapping
$(\eps,\theta)\mapsto\rho_\eps(\theta)$. The first-order expansion
in \eqref{fin} and Taylor's theorem give
\begin{equation}\label{eq:rho-initial}
\rho_\eps=1+\eps\cos 2\theta+\eps^2h_2+O(\eps^3),
\end{equation}
with a remainder uniform in $\theta$. The wave numbers are real
analytic and, by \eqref{eveness}, even in $\eps$. Thus we also have
\begin{equation}\label{eq:frequency-quadratic}
k_\eps=k_*+k^{(2)}\eps^2+O(\eps^4),
\qquad
\kappa_\eps=\kappa_*+\kappa^{(2)}\eps^2+O(\eps^4).
\end{equation}
The reflection symmetries make $h_2$ even and $\pi$-periodic.
The exact normalisation
\eqref{eq:polar-normalisation} gives
\begin{equation}\label{eq:normalisation-consequences}
\langle\rho_\eps-1\rangle=0,
\qquad
\langle(\rho_\eps-1)\cos 2\theta\rangle=\frac{\eps}{2},
\qquad
\langle h_2\rangle=\langle h_2\cos 2\theta\rangle=0.
\end{equation}
Taylor expansion of $F_n$ in $\rho$ about $\rho=1$, at
$(k,\kappa)=(k_\eps,\kappa_\eps)$, gives
\begin{align}
M_n
&=
\delta_{n0}F_0(1;k_\eps,\kappa_\eps)
+b_n(1;k_\eps,\kappa_\eps)
\langle(\rho_\eps-1)\cos(n\theta)\rangle
\nonumber\\
&\quad+
\frac{\partial_rb_n(1;k_\eps,\kappa_\eps)}{2}
\langle(\rho_\eps-1)^2\cos(n\theta)\rangle
\label{eq:radial-taylor}\\
&\quad+
\frac{\partial_r^2b_n(1;k_\eps,\kappa_\eps)}{6}
\langle(\rho_\eps-1)^3\cos(n\theta)\rangle
+O(\eps^4),
\nonumber
\end{align}
where $\delta_{n0}=1$ for $n=0$ and $0$ otherwise. For $n>0$,
the term $F_n(1;k_\eps,\kappa_\eps)$ is independent of $\theta$
and therefore has zero average against $\cos(n\theta)$.
Below, $b_n$, $b_n'$ and $b_n''$ denote $b_n$ and its first two
radial derivatives at $(r,k,\kappa)=(1,k_*,\kappa_*)$.
The properties of $b_n(r;k,\kappa)$ that we need are given in
Lemma~\ref{lem:bessel-data}. For $j\ge2$, the
coefficient of $\eps^2$ in $M_{2j}=0$ is
\[
b_{2j}\langle h_2\cos(2j\theta)\rangle
+\frac{b_{2j}'}{2}\langle(\cos 2\theta)^2\cos(2j\theta)\rangle=0.
\]
For $j\ge3$, the second average vanishes and $b_{2j}\ne0$.
Together with \eqref{eq:normalisation-consequences} and the
symmetries of $h_2$, this gives $h_2=A\cos 4\theta$ for some
$A\in\R$. Taking $j=2$ gives $Ab_4/2+b_4'/8=0$, so
Lemma~\ref{lem:bessel-data}(c) yields
\begin{equation}\label{eq:shape-coefficient}
A=-\frac{b_4'}{4b_4}=\frac34+\frac{\kappa_*^2}{12}.
\end{equation}
Substituting this value into \eqref{eq:rho-initial} gives
\eqref{eq:mainshape}.

We now turn to the expansion \eqref{eq:main-contrast} of the material coefficient. Set $d_*:=d_0(k_*)=d_0(\kappa_*)=-\kappa_*^2/2$, see \eqref{eq:seed-d0}, and recall that $\kappa^{(2)}$ and $k^{(2)}$ are the coefficients in \eqref{eq:frequency-quadratic}. The cubic term of $M_2=0$ determines $\kappa^{(2)}$. Since $J_2(\kappa_*)=0$, \eqref{eq:frequency-quadratic} and Lemma~\ref{lem:bessel-data}(a) give
\[
b_2(1;k_\eps,\kappa_\eps)
=J_0(k_\eps)J_2(\kappa_\eps)
=\frac{\kappa^{(2)}}{\kappa_*}b_2'\eps^2+O(\eps^4).
\]
By \eqref{eq:normalisation-consequences}, $\langle(\rho_\eps-1)\cos 2\theta\rangle=\eps/2$. Using $h_2=A\cos 4\theta$, we have
\[
\langle(\rho_\eps-1)^2\cos 2\theta\rangle
=\frac{A}{2}\eps^3+O(\eps^4),
\qquad
\langle(\rho_\eps-1)^3\cos 2\theta\rangle
=\frac38\eps^3+O(\eps^4).
\]
Thus the coefficient of $\eps^3$ in \eqref{eq:radial-taylor} for $n=2$ gives
\begin{equation}\label{eq:M2-cubic}
\frac{\kappa^{(2)}}{2\kappa_*}b_2'
+\frac A4b_2'
+\frac1{16}b_2''=0.
\end{equation}
Dividing by $b_2'/2$ and using Lemma~\ref{lem:bessel-data}(a) and \eqref{eq:shape-coefficient}, we obtain
\begin{equation}\label{eq:kappa2}
\frac{\kappa^{(2)}}{\kappa_*}
=-\frac A2-\frac{1+2d_*}{8}
=\frac{\kappa_*^2-6}{12}.
\end{equation}
To determine $k^{(2)}$, take the coefficient of $\eps^2$ in $M_0=0$. Since $\langle\rho_\eps-1\rangle=0$ and $\langle(\rho_\eps-1)^2\rangle=\eps^2/2+O(\eps^3)$, \eqref{eq:radial-wronskian} and Lemma~\ref{lem:bessel-data}(b) give
\begin{equation}\label{eq:k2-equation}
d_0'(k_*)k^{(2)}-d_0'(\kappa_*)\kappa^{(2)}
=-\frac{(\kappa_*^2-k_*^2)(1+2d_*)}{4}
=\frac{(\kappa_*^2-k_*^2)(\kappa_*^2-1)}{4}.
\end{equation}
From \eqref{eq:d0-derivative} and \eqref{eq:seed-d0} we compute
\[
d_0'(\kappa_*)=-\frac{\kappa_*(\kappa_*^2+4)}{4},
\qquad
d_0'(k_*)=-\frac{\kappa_*^4+4k_*^2}{4k_*}.
\]
Substituting these values and \eqref{eq:kappa2} into
\eqref{eq:k2-equation} yields
\begin{equation}\label{eq:k2}
\frac{k^{(2)}}{k_*}
=
\frac{
\kappa_*^2(\kappa_*^2-6)(\kappa_*^2+4)/12
-(\kappa_*^2-k_*^2)(\kappa_*^2-1)
}{
\kappa_*^4+4k_*^2
}.
\end{equation}
Finally, \eqref{eq:frequency-quadratic} gives
\begin{equation}\label{eq:contrast-expansion}
\begin{split}
\sigma_\eps
&=\frac{\kappa_\eps^2}{k_\eps^2}
=\sigma_*+c_*\eps^2+O(\eps^4),\\
c_*&=2\sigma_*
\left(\frac{\kappa^{(2)}}{\kappa_*}-\frac{k^{(2)}}{k_*}\right)
=\frac{2\sigma_*(2\kappa_*^2+3)(\kappa_*^2-k_*^2)}
{3(\kappa_*^4+4k_*^2)}.
\end{split}
\end{equation}
Since $\kappa_*>k_*$, we have $c_*>0$.

It remains to show that the domains $\Omega_\eps$ are not ellipses.
\begin{proposition}\label{prop:nonellipticity}
For every sufficiently small nonzero $\eps$, the domain $\Omega_\eps$
is not an ellipse.
\end{proposition}

\begin{proof}
The reflection symmetries would force an ellipse in the family to
be centred at the origin with axes along the coordinate axes.
Its polar radius would therefore be $C(1-\eta\cos 2\theta)^{-1/2}$,
with $C>0$ and $\eta$ small, since $\rho_\eps\to1$ uniformly. From
\[
(1-\eta\cos 2\theta)^{-1/2}
=
1+\frac{\eta}{2}\cos 2\theta
+\frac{3\eta^2}{16}(1+\cos 4\theta)+O(\eta^3),
\]
the normalisations \eqref{eq:polar-normalisation} give
$C=1-3\eta^2/16+O(\eta^3)$ and $\eps=\eta/2+O(\eta^3)$.
The polar radius of the normalised ellipse would thus be
\[
\rho_{\mathrm{ell}}(\theta)
=
1+\eps\cos 2\theta+\frac34\eps^2\cos 4\theta+O(\eps^3).
\]
By \eqref{eq:mainshape}, the $\cos 4\theta$ coefficients of
$\rho_\eps$ and $\rho_{\mathrm{ell}}$ differ by
$\kappa_*^2\eps^2/12+O(\eps^3)$, which is nonzero for sufficiently
small $\eps\ne0$.
\end{proof}

\begin{proof}[Proof of Theorem~\ref{thm:non-scattering-branch}]
After reducing $\eps_0$ if necessary, the normalised family \eqref{eq:normalized-family} has all the
stated properties by Propositions~\ref{prop:bifurcating-branch},
\ref{prop:branch-geometry} and~\ref{prop:nonellipticity}, together
with dilation invariance.
\end{proof}

\begin{proof}[Proof of Proposition~\ref{prop:normalised-expansions}]
The normalisation was established in \S\ref{sec:geometry}.
The shape expansion follows from \eqref{eq:rho-initial} and
\eqref{eq:shape-coefficient}, and the material coefficient expansion is
\eqref{eq:contrast-expansion}.
\end{proof}

%%%%%%%%%%%%%%%%%%%%%%%%%%%%%
\section*{\phantomsection Acknowledgements}
\addcontentsline{toc}{section}{Acknowledgements}
The work of FC was partially supported by NSF grant DMS-2406313. The work of ML was partially supported by the AI for Math Fund, an initiative by Renaissance Philanthropy funded by XTX Markets. 
The authors would like to thank the Isaac Newton Institute for Mathematical Sciences, Cambridge, where
this investigation started, for support and hospitality during the programme \emph{Geometric spectral
theory and applications}, supported by EPSRC grant EP/Z000580/1. ML's stay at INI was supported by a
Simons fellowship.

%%%%%%%%%%%%%%%%%%%%%%%%%%%%%
\section*{\phantomsection AI usage statement}
\addcontentsline{toc}{section}{AI usage statement}
The main idea of the paper arose from discussions with ChatGPT 6 and Claude Fable 5.2 about the possibility of adapting the constructions of Schiffer domains to obtain examples of non-scattering domains.  
Claude was also used for proofreading and to produce a \texttt{Python} script mentioned below. The authors take full responsibility for the contents and presentation of the paper.

%%%%%%%%%%%%%%%%%%%%%%%%%%%%%
\section*{\phantomsection Data availability statement}
\addcontentsline{toc}{section}{Data availability statement}
The accompanying  \texttt{Python} and \texttt{Mathematica} scripts used to produce Figure~\ref{fig:nonscattering-domain} are available for download at
\begin{center}
{\small
\url{https://www.michaellevitin.net/non-scattering.html}\quad and\quad \url{https://github.com/michaellevitin/non-scattering}.}
\end{center}

%%%%%%%%%%%%%%%%%%%%%%%%%%%%%
\appendix
\section{Auxiliary results}
We collect here the auxiliary results used in the paper.
%%%%%%%%%%%%%%%%%%%%%%%%
\subsection{The analytic form of the Crandall--Rabinowitz theorem}\label{A2}
\begin{theorem}
\label{thm:analytic-cr}
Let $X$ and $Y$ be real Banach spaces, let $U\subset X$ be an open
neighbourhood of $0$, and let $I\subset\R$ be an open interval
containing $\lambda_*$. Suppose that $F:U\times I\to Y$ is real
analytic and that $F(0,\lambda)=0$ for every $\lambda\in I$.
Assume that $L=D_xF(0,\lambda_*)$ has kernel
$\operatorname{span}\{\phi\}$, with $\phi\ne0$, and closed range
of codimension one. If
\[
D_\lambda D_xF(0,\lambda_*)\phi\notin\Ran L,
\]
then there exist $\delta>0$ and a real-analytic curve
$(-\delta,\delta)\ni s\mapsto(x(s),\lambda(s))\in U\times I$
such that
\[
F(x(s),\lambda(s))=0,
\qquad
\lambda(0)=\lambda_*,
\qquad
x(s)=s\phi+O_X(s^2).
\]
Moreover, every zero of $F$ sufficiently close to $(0,\lambda_*)$
either has $x=0$ or lies on this curve.
\end{theorem}
For the proof we refer the reader to \cite[Theorem~8.3.1]{BT}.
%%%%%%%%%%%%%%%%%%%%%%%%%%%%%%%%%%%%%%%%%%%%%%%
\subsection{Bessel identities}
Recall that $b_n$, $b_n'$ and $b_n''$ denote the values of $b_n(r;k,\kappa)$, see \eqref{eq:radial-integrals}, and of its first two radial derivatives at $(r,k,\kappa)=(1,k_*,\kappa_*)$, and that $d_*=-\kappa_*^2/2$.
\begin{lemma}\label{lem:bessel-data}
We have
\begin{enumerate}[\rm(a)]
\item $b_2=0$, $b_2'=\kappa_*J_0(k_*)J_2'(\kappa_*)\ne0$, and
$b_2''=(1+2d_*)b_2'$.
\item $b_0'=J_0(k_*)J_0(\kappa_*)(1+2d_*)$.
\item $b_{2j}\ne0$ for $j\ge2$, and
$b_4'/b_4=-3-\kappa_*^2/3$.
\item For $k\ne\kappa$ with $J_0(k)J_0(\kappa)\ne0$,
\begin{equation}\label{eq:radial-wronskian}
F_0(1;k,\kappa)
=
\frac{J_0(k)J_0(\kappa)}{\kappa^2-k^2}
\left(d_0(k)-d_0(\kappa)\right).
\end{equation}
\end{enumerate}
\end{lemma}

\begin{proof}
(a) Since $J_2(\kappa_*)=0$, we have $b_2=0$. Differentiating
$rJ_0(k_*r)J_2(\kappa_*r)$ at $r=1$ gives the stated value of
$b_2'$, which is nonzero because $J_0(k_*)\ne0$ and $\kappa_*$ is
a simple zero of $J_2$. Differentiating again and using
$\kappa_*^2J_2''(\kappa_*)=-\kappa_*J_2'(\kappa_*)$, by Bessel's
equation, gives
\[
b_2''
=\kappa_*J_2'(\kappa_*)
\left(J_0(k_*)+2k_*J_0'(k_*)\right)
=(1+2d_*)b_2'.
\]

(b) Direct differentiation gives
\[
b_0'
=J_0(k_*)J_0(\kappa_*)
\left(1+d_0(k_*)+d_0(\kappa_*)\right).
\]

(c) For $j\ge2$, $J_{2j}(\kappa_*)\ne0$ by the argument in the
proof of Proposition~\ref{prop:kernel-range}. Since $J_0(k_*)\ne0$,
we have $b_{2j}\ne0$. The recurrences
$J_{n+1}(t)=2nJ_n(t)/t-J_{n-1}(t)$ and
$J_n'(t)=J_{n-1}(t)-nJ_n(t)/t$, together with
$J_2(\kappa_*)=0$ and
$J_1(\kappa_*)=\kappa_*J_0(\kappa_*)/2$, give
\[
J_3(\kappa_*)=-\frac{\kappa_*}{2}J_0(\kappa_*),
\qquad
J_4(\kappa_*)=-3J_0(\kappa_*),
\qquad
\frac{\kappa_*J_4'(\kappa_*)}{J_4(\kappa_*)}
=\frac{\kappa_*^2}{6}-4.
\]
Therefore
\[
\frac{b_4'}{b_4}
=1+d_*+\frac{\kappa_*J_4'(\kappa_*)}{J_4(\kappa_*)}
=-3-\frac{\kappa_*^2}{3}.
\]

(d) Integrating the radial Wronskian identity for $J_0(kr)$ and
$J_0(\kappa r)$ gives
\[
(\kappa^2-k^2)F_0(1;k,\kappa)
=
kJ_0'(k)J_0(\kappa)-\kappa J_0(k)J_0'(\kappa).
\]
The definition of $d_0$ now gives \eqref{eq:radial-wronskian}.
\end{proof}

\end{document}